\documentclass[11pt, oneside]{article}
\usepackage{amsfonts}
\usepackage{mathrsfs}
 \usepackage{color}
\usepackage{latexsym}
\usepackage{amssymb}
\usepackage{amsmath}
\usepackage{enumerate}
\usepackage{amsthm}
\usepackage{indentfirst}
\usepackage{mathtools}
\usepackage{tikz}
\usepackage{psfrag}
\usepackage{graphicx}
\usepackage{bm}
\usepackage[rflt]{floatflt}
\usepackage{float}
\usepackage{authblk}
\usepackage{ulem}
\usepackage[square, comma, sort&compress, numbers]{natbib}
\usepackage{verbatim}
\usepackage{amsthm}
\usepackage{esint}

\numberwithin{equation}{section}
\allowdisplaybreaks

\newenvironment{proof2.1}{\medskip\noindent{\bf Proof of the Theorem 2.1:}\enspace}{\hfill \qed \newline \medskip}

\newenvironment{proof2.2}{\medskip\noindent{\bf Proof of the Theorem 2.2:}\enspace}{\hfill \qed \newline \medskip}

\newtheorem{theorem}{\color{black}\indent Theorem}[section]
\newtheorem{lemma}{\color{black}\indent Lemma}[section]

\newtheorem{definition}{\color{black}\indent Definition}[section]
\newtheorem{remark}{\color{black}\indent Remark}[section]

\usepackage{amssymb,amsmath}
\begin{document}
\title{On higher integrability  for  $p(x)$-Laplacian equations with   drift}
\author[1]{Jingya Chen}
\author[1]{Bin Guo\thanks{Corresponding author
\newline \hspace*{4mm}{
Email addresses: bguo@jlu.edu.cn (Bin Guo)}}}
\author[2]{Baisheng Yan}
\affil[1]{School of Mathematics, Jilin University, Changchun, Jilin Province 130012, China}
\affil[2]{Department of Mathematics, Michigan State University, East Lansing,
MI 48824, USA}
\renewcommand*{\Affilfont}{\small\it}
\date{} \maketitle
\vspace{-20pt}

{\bf Abstract:}  In this paper, we study the higher integrability for the gradient  of weak solutions of  $p(x)$-Laplacians equation with drift terms.  We prove a version  of generalized  Gehring's lemma under some weaker condition on  the modulus of continuity of variable exponent $p(x)$ and present a modified version of Sobolev-Poincar\'{e} inequality with such an exponent. When $p(x)>2$ we derive the reverse H\"older inequality with a proper dependence on the drift and force terms and establish a specific high  integrability result. Our condition on the exponent $p(x)$ is more specific and weaker than the known conditions and our results extend some  results  on the $p(x)$-Laplacian equations  without  drift terms. \\


\section{Introduction}
The classical $p$-Laplacian equation $\mbox{div} (|\nabla u|^{p-2}\nabla u)=0,$ where $p>1$ is a constant, has been widely studied, including   some estimates for very weak solutions of $p$-Laplacian equations below the natural energy space $W^{1,p}.$ We refer to \cite{byu1,caff,ct,dib,iw1,is, kinn2, le, ya}   for some  results and further references. As a natural generalization of the $p$-Laplacian equation intended for many important applications,  the $p(x)$-Laplacian equation
\begin{align}\label{no-dri}\hbox{div}(|\nabla u|^{p(x)-2}\nabla u)=0, \quad p(x)>1,
\end{align}
provides  good models for certain nonlinear and nonhomogenous features of some materials, such as  the electro rheological fluids  \cite{rr1,rr2},  elastic mechanics \cite{Zhi2}, image restoration \cite{clr} and flows in porous media \cite{as1}. There have been extensive and  growing literatures associated with the $p(x)$-Laplacian equation; for some  overview, results and further references, we refer to \cite{ace,yay,byu1,crz,die2,die1, fmg,gmvr,Zhi2,Zhi1}. For example, qualitative theories for  $p(x)$-Laplacian  equations  have been  developed under the so-called ``logarithmic condition" on the exponent $p(x)$;  i.e., for some $k>0,$
\begin{align}\label{log}
|p(x)-p(y)|\leqslant \frac{k}{\ln\frac{1}{|x-y|}}, \quad 0<|x-y|<1.
\end{align}
 It is well known that the condition \eqref{log} plays an important role in the study of  the Hardy-Littlewood maximal function and the Riesz integral operators in the relevant Lebesgue-Orlicz space \cite{die2}, the higher integrability of  gradient \cite{die1}, and the H\"{o}lder continuity  of $p(x)$-harmonic functions \cite{phph1,phph2}.

  In this paper, we are  mainly concerned with the higher integrability for the gradient of weak solutions to the  Dirichlet problem:
\begin{equation}\label{equ}
\begin{cases}\text{div}(|\nabla  u|^{p(x)-2}\nabla  u)+\vec{b}(x)\cdot\nabla  u=\hbox{div}(\vec{h}) \quad \text{in}~\Omega,\\ u(x)=0,~x\in\partial\Omega,
\end{cases}\end{equation}
where $\Omega\subset\mathbb{R}^n$  is a bounded domain  satisfying the  \textit{uniform exterior cone condition} (see Definition \ref{cone-Omega} below),  and $\vec{b}(x)$ and $\vec{h}(x)$ are given $n$-dimensional vector functions. By a  weak  solution to problem  (\ref{equ})~we mean a function
\[
u\in W^{1,p(\cdot)}_0(\Omega)=\{u\in W^{1,1}_0(\Omega): \int_\Omega|\nabla  u|^{p(x)}dx<\infty\}
\]
 such that the integral identity
\begin{equation}\label{defweak}
\int_{\Omega}\Big [|\nabla  u|^{p(x)-2}\nabla  u\cdot\nabla \varphi-\varphi \vec{b}(x)\cdot\nabla  u
 \Big]dx=\int_{\Omega}\vec{h}\cdot\nabla \varphi dx
\end{equation}
holds for all test functions  $\varphi\in W^{1,p(\cdot)}_0(\Omega).$ We refer to \cite{crz,die2,Zhi1} for more information on the general Lebesgue-Orlicz spaces including the Sobolev-Orlicz space $W^{1,p(x)}_0(\Omega).$

The study of higher integrability for the gradient of weak solutions to (\ref{equ})  relies  essentially   on establishing a modified version of   reverse H\"older inequality of the form:
\begin{align}\label{addrh}
\fint_{B_{r}}|f|^{q}dx\leqslant C_1\Big(\fint_{B_{2r}}|f|dx\Big)^{q+\omega(r)}+C_2\fint_{B_{2r}}Gdx \qquad\mbox{for all balls $B_r$,}
\end{align}
for $f=|\nabla u|^{\frac{p(x)}{q}}$ with a constant $q>1$,   a continuous increasing function $\omega(r)$  with $\omega(0)=0$ and a nonnegative function $G\in L^1(\Omega).$ Such an inequality is usually  derived from  the standard  Cacciopolli's inequalities and certain modified version of Sobolev-Poinca\'re inequalities.  Under the logarithmic condition  \eqref{log},   function $\omega(r)$ can be chosen such that  $r^{-n\omega(r)}\leqslant C'.$ In this way,  \eqref{addrh} becomes the classical reverse H\"older inequality, from which the higher integrability follows by the well-known Gehring lemma \cite{geh}.   Without the log-continuity, the space $W^{1,p(x)}_0(\Omega)$ may not behave like the usual Sobolev space $W^{1,p}_0(\Omega)$ with constant $p$ and many traditional approaches for  elliptic problems do not work anymore.  Nevertheless, a nonlogarithmic condition has been studied by Zhikov \cite{Zhi2} to  guarantee  the absence of Lavrentiev's phenomenon. In particular, such a result holds if $p(x)$ satisfies
 \begin{align}\label{nolog}
 |p(x)-p(y)|\leqslant  L\frac{\log\log\frac{1}{|x-y|}}{\log\frac{1}{|x-y|}}, \quad 0<L\leqslant \frac{1}{n},~~0<|x-y|\ll1.
 \end{align}
Later, under this condition Zhikov and Pastukhova   in \cite{Zhi1} established  a  version of generalized  Gehring's lemma and then proved a higher integrability for the gradient of  weak solutions of  \eqref{equ} without the drift term $\vec{b}.$

We revisit the work of Zhikov and Pastukhova  \cite{Zhi1} and prove a generalized  Gehring's lemma under a weaker and more specific condition on  the modulus of variable exponent $p(x)$ and present a modified version of Sobolev-Poincar\'{e} inequality with such an exponent. Finally we derive the reverse H\"older inequality with a proper dependence on the drift and force terms, which enables us to  establish a  specific high  integrability result. Our   result  extends the corresponding  result of \cite{Zhi1}  on the $p(x)$-Laplacian equations  without  drift.
Throughout  the paper,  we assume   that there exist numbers $\beta\geqslant \alpha>1$ and $r_0 \in (0,1)$  and a  non-decreasing  continuous function $\omega:[0,+\infty)\to[0,+\infty)$  such that
\begin{equation}\label{omeg}
\begin{aligned}
& (a) \quad  \alpha\leqslant p(x)\leqslant\beta \quad \forall\, x\in\Omega,\\
&  (b) \quad   |p(x)-p(y)|\leqslant\omega\big ( {|x-y|}/{2}\big ) \quad  \forall\, x,y\in\Omega,\\
&  (c) \quad  \omega'(r)\geqslant  0 \quad \forall \, r\in (0,r_0],\\
&  (d) \quad  \omega(0)=0, \quad \omega_0=\omega(r_0)<1.
\end{aligned}
\end{equation}

To state our main assumption and result, given $q>1$, let $ \lambda_0=r_0^{-n/q}>1$  and define $a(\lambda)= \lambda^{\omega(\lambda^{-\frac{q}{n}})} $ for all $\lambda\geqslant \lambda_0.$  Then    it can be easily shown that the function $s= \frac{\lambda}{a(\lambda)}$   is invertible  from   $ [\lambda_0,\infty)$ onto $[s_0,\infty)$  and let  $\lambda=\lambda(s)$ be its inverse function on $[s_0,\infty).$
Our specific assumption on the modulus of continuity $\omega$ of the exponent $p(x)$ is the following:
\begin{equation}\label{ass-0}
 \sup \left\{ \int_{s}^{\lambda(ks)} \frac{1}{t\,a^q(t) } dt :s\geqslant \lambda_0 \right\} <\infty \quad \quad \forall\, k\geqslant 1.
\end{equation}
Furthermore, given  $\theta>0$, let  $\Phi=\Phi_{\theta}\colon [0,\infty)\to (0,\infty)$ be a  $C^1$ increasing function such that
\begin{equation}\label{Phi}
\Phi(\lambda)= \exp\left(\theta (1- \omega_0)\int_{\lambda_0}^{\lambda}\frac{1}{ta^q(t)}dt\right) \qquad \forall\, \lambda \geqslant \lambda_0.
\end{equation}
The main result of the paper is the following theorem.

\begin{theorem}\label{main3} Suppose that
 assumptions   (\ref{omeg}) and (\ref{ass-0}) are satisfied for some $\alpha>2.$ Let \[
G=1+|\vec{b}(x)|^{\frac{p(x)}{p(x)-2}}+ |\vec{h}(x)|^{\frac{p(x)}{p(x)-1}}\in L^1(\Omega).
\]
Then there exist  numbers $q>1$  and $\theta>0$ such that for a function $\Phi=\Phi_\theta$ as  given above if
\begin{equation}\label{g-i}
G(x)\Phi(G(x)^\frac{1}{q})\in L^1(\Omega),
\end{equation}
then any weak solution $u$ to the problem (\ref{equ}) will satisfy
\begin{equation}\label{h-i}
|\nabla  u|^{p(x)}\Phi(|\nabla  u|^\frac{p(x)}{q})\in L^1(\Omega).
\end{equation}
\end{theorem}

\begin{remark}  Note that  condition (\ref{g-i}) is automatically satisfied if both $\vec{b}$ and $\vec{h}$ belong to $L^\infty(\Omega).$ Moreover, the conclusion (\ref{h-i}) indeed guarantees  a higher integrability for the gradient $\nabla u$ provided that
\begin{equation}\label{ass-11}
\int_{\lambda_0}^{\infty}\frac{1}{ta^q(t)}dt =\infty.
\end{equation}
 Furthermore, it can be verified that  the assumptions (\ref{ass-0}) and (\ref{ass-11})  both hold for the function
\[
\omega(r)=\frac{k}{\ln \frac{1}{r}} \;\; \mbox{ with }\; k\geqslant 0,  \quad \mbox{ or}\quad \omega(r)=L\frac{\ln\ln\frac{1}{r}}{\ln\frac{1}{r}} \;\; \mbox{ with } \;  0<L\leqslant \frac{1}{n}.
\]
\end{remark}

The outline of this paper is the following: In Section~2, we examine some consequences of the condition (\ref{ass-0}) (Lemma \ref{lem-21}) and  prove a  version of generalized  Gehring's lemma (Theorem \ref{gehring0}) that improves the result of \cite[Theorem 6.4]{Zhi1}. In Section 3, we prove  a modified version Sobolev-Poincar\'{e} inequality  involving the modulus of  exponent $p(x)$ (Theorem \ref{gesobpoi1}). In Section 4, we establish the Cacciopolli inequality with specific dependence on the drift and force terms (Theorem \ref{cacc}) and for  domains satisfying the uniform exterior cone condition we derive  the  reverse H\"older inequality  (Theorem \ref{reHolder2}), from which  our main theorem follows.

\section{Generalized Gehring's  lemma}

In what follows, we assume that  $\omega:[0,+\infty)\to[0,+\infty)$   is  a non-decreasing  continuous function which is also differentiable on $(0,r_0]$ for some $r_0\in (0,1)$. Moreover, we  assume that $\omega$ satisfies the conditions (c) and (d) of  (\ref{omeg}).
For $q>1$, we  extend the functions $a(\lambda)$ and $s(\lambda)$ defined above to all $\lambda\geqslant 0$ by setting \begin{equation}\label{adef0}
a(\lambda)=\begin{cases} \lambda^{\omega(\lambda^{-\frac{q}{n}})} &\forall \, \lambda\in [\lambda_0,\infty) \\ \lambda_0^{\omega_0} &\forall\,\lambda\in [0,\lambda_0) \end{cases} \quad\mbox{and}   \quad s(\lambda)=\frac{\lambda}{a(\lambda)}, \quad \mbox{for all $ \lambda \geqslant  0.$ }
\end{equation}
Note that if $\lambda\geqslant \lambda_0$ then   $
s(\lambda)=\lambda^{1-\omega(\lambda^{-\frac{q}{n}})}\geqslant \lambda^{1-\omega_0}$ and
\[
s^\prime(\lambda) =\frac{1-\omega(\lambda^{-\frac{q}{n}})+\frac{q}{n}\ln\lambda\cdot \omega^\prime(\lambda^{-\frac{q}{n}})\lambda^{-\frac{q}{n}}}{ {a}(\lambda)}\geqslant \frac{1-\omega_0}{a(\lambda)} >0 .
\]
This shows that the function $s=s(\lambda)$  is invertible  from  $ [\lambda_0,\infty)$ onto $[s_0,\infty),$ where $s_0=\lambda_0^{1-\omega_0}.$ Let   $\lambda=\lambda(s)$ be the inverse function of $s=s(\lambda)$,  from  $[s_0,\infty)$ onto $ [\lambda_0,\infty)$.  Then
\begin{equation}\label{d-lambda}
0< \lambda^\prime(s)=\frac{1}{s^\prime(\lambda)}\leqslant\frac{ {a}(\lambda(s))}{1-\omega_0} \qquad \forall\,s\geqslant s_0.
\end{equation}

\begin{lemma}\label{lem-21} Let $q>1$ and assume the condition (\ref{ass-0}) holds. Then for all $k>0$,
\begin{equation}\label{ass-1}
M(k):=\sup  \left\{ \int_{s}^{\lambda(ks)} \frac{1}{t\,a^q(t) } dt : s\geqslant \frac{s_0}{k}\right\} <\infty.
\end{equation}
Given  $\theta>0$ and $c>0$, let $\varphi=\varphi_{\theta,c}\colon [0,\infty)\to (0,\infty)$ be an increasing $C^1$ function such that
\begin{equation}\label{phidef0}
 {\varphi}  (s)=
\exp\left(\theta (1-\omega_0)\int_{\frac{s_0}{4c}}^{s}\frac{1}{ t{a}^q(t) } dt \right) \quad\mbox{for all } \; s\geqslant \frac{s_0}{4c}.  \end{equation}
Then for each $k>0,$ there exists a constant $L(k)>0$ such that
\begin{equation}\label{step4}
 \varphi(ks)\leqslant L(k)\varphi(s)\quad\quad\forall\,s\geqslant 0.
\end{equation}
Moreover, if $\Phi=\Phi_{\theta}$ is a function  as given in (\ref{Phi}), then there exist  $\nu_2>\nu_1>0$ such that
\begin{equation}\label{step2}
\nu_1\leqslant  {\Phi(s)}/{\varphi(s)}\leqslant \nu_2 \quad\quad\forall\,s\geqslant 0.
\end{equation}
\end{lemma}

\begin{proof}  If $k\geqslant 1,$ then for all $s\geqslant \frac{s_0}{k}$, by   (\ref{ass-0}),
\[
 \int_{s}^{\lambda(ks)}  \frac{dt}{t\,a^q(t) } \leqslant    \max\left\{ \sup_{s\geqslant \lambda_0}  \Big (\int_{s}^{\lambda(ks)} \frac{dt}{t\,a^q(t) } \Big), \; \max_{\frac{s_0}{k}\leqslant s \leqslant \lambda_0} \Big (\int_{s}^{\lambda(ks)} \frac{dt}{t\,a^q(t) }  \Big ) \right\} <\infty;
\]
 if $0<k<1,$ then for all $s\geqslant \frac{s_0}{k}$, by   (\ref{ass-0}) with $k=1,$
\[
\begin{split} &\int_{s}^{\lambda(ks)}  \frac{dt}{t\,a^q(t) } \leqslant  \int_{ks}^{\lambda(ks)} \frac{dt}{t\,a^q(t) }  \leqslant  \max\left\{ \sup_{s'\geqslant \lambda_0}  \Big (\int_{s'}^{\lambda(s')} \frac{dt}{t\,a^q(t) } \Big), \; \max_{s_0\leqslant s'\leqslant \lambda_0} \Big (\int_{s'}^{\lambda(s')} \frac{dt}{t\,a^q(t) }  \Big ) \right\} <\infty.
 \end{split}
\]
This proves  (\ref{ass-1}). Next, clearly (\ref{step4}) holds with $L(k)=1$ if $0<k\leqslant 1,$ so we assume $k>1$ and   prove (\ref{step4})
 in two cases:   (i) $1< k\leqslant 4c,$ and  (ii) $k>4c .$
In Case (i),
since  $\varphi$ is increasing and positive, by virtue  of  $\lambda(s)\geqslant s$ for all $s\geqslant s_0$, we have
$
\varphi(ks)\leqslant \varphi(4c s) \leqslant \varphi(\lambda(4c s))\le e^{\theta M(4c )}\varphi(s)$ for all $s\geqslant \frac{s_0}{4c}.
$
If $0\leqslant s\leqslant \frac{s_0}{4c}<s_0,$ then $\varphi(ks)\leqslant \varphi( {ks_0} )\leqslant  {\varphi(ks_0)}\varphi(s)/{\varphi(0)}.$ Thus   (\ref{step4}) follows in this case with
\[
L(k)=\max\Big\{ {\varphi( ks_0)}/{\varphi(0)}, \; e^{\theta M(4c)}\Big\}.
\]
In Case (ii), if $s\geqslant \frac{s_0}{4c},$ then $s> \frac{s_0}{k}$  and thus  $\varphi(ks)\leqslant  \varphi(\lambda(ks))\le e^{\theta M(k)}\varphi(s).$ If $0\leqslant s\leqslant \frac{s_0}{4c},$ then $\varphi(ks)\leqslant \varphi(\frac{ks_0}{4c})\leqslant  {\varphi(\frac{ks_0}{4c})}\varphi(s)/{\varphi(0)}.$  Thus   (\ref{step4}) follows in this case with
\[
L(k)= \max\Big\{ {\varphi(\frac{ks_0}{4c})}/{\varphi(0)}, \,e^{\theta M(k)}\Big\}.
\]
This proves (\ref{step4}). To prove (\ref{step2}), we note that
\[
\frac{\Phi(s)}{\varphi(s)} =\exp\left (\theta(1-\omega_0) \int_{\lambda_0}^{\frac{s_0}{4c}} \frac{dt}{ta^q(t)}\right) := \bar\nu >0 \quad\forall\, s\geqslant \bar s:=\max\{\lambda_0,\frac{s_0}{4c}\}.
\]
Thus (\ref{step2}) follows with
\[
\nu_1=\min\left \{\bar\nu,\;  \min_{s\in [0,\bar s]} \frac{\Phi(s)}{\varphi(s)} \right \},\qquad \nu_2=\max\left \{\bar\nu, \; \max_{s\in [0,\bar s]} \frac{\Phi(s)}{\varphi(s)} \right \}.
\]
This completes the proof of the lemma.
\end{proof}

  In the following, we denote by  $Q_r= Q_r(x_0)=x_0+[-r/2,r/2]^n$ the cube with edge $r>0$~and  center  $x_0\in\mathbb{R}^n.$ Note that $|Q_r|=r^n.$

\begin{theorem}\label{gehring0}  Let $q>1$ and  let $f$ and $g\in L^q(\mathbb{R}^n)$~be non-negative functions with
\begin{equation}\label{norm-1}
\|f\|_q\leqslant 1
\end{equation}
such that for some constants $l\geqslant 1$, $c_0>0$ and $c_1>0,$  the inequality
\begin{equation}\label{gecod0}
\fint_{Q_r(x_0)}f^qdx\leqslant c_0\left(\fint_{Q_{lr}(x_0)}fdx\right)^{q+\omega(r)}+c_1\fint_{Q_{lr}(x_0)}g^qdx
\end{equation}
holds for all $ 0<r\leqslant r_0$ and $x_0\in\mathbb{R}^n.$  Then it follows that
\begin{equation}\label{result0}
\int_{\mathbb{R}^n}f^q {\varphi}(f)dx\leqslant 4 {\varphi}(\lambda_0)+C\int_{\mathbb{R}^n}g^q {\varphi}(g)dx
\end{equation}
for a constant $C>0$ independent of $l, f$ and $g$, where  $\varphi=\varphi_{\theta,c_0}\colon [0,\infty)\to (0,\infty)$ is an increasing  $C^1$ function of the form (\ref{phidef0}) with $\theta>0$ being the constant determined by
\begin{equation}\label{theco0}
 \frac{4^{n}(4c_0)^{q}\theta e^{\theta M(4c_0)} }{q-1}=\frac{1}{2}.
\end{equation}
\end{theorem}
\begin{proof} We prove the theorem by following  the main idea of~\cite[Theorem~6.4]{Zhi1}  with some new observations. To do this, we split the proof into several steps.\\

{\bf Step 1.}   Let $\lambda \geqslant \lambda_0.$ We claim that for almost all~$x_0\in E(\lambda)$ there exists~$r=r(x_0)\in (0,r_0]$~such that
\begin{equation}\label{e}
\fint_{Q_{s}}f^qdx\leqslant\fint_{Q_r}f^qdx=\lambda^q\qquad \forall\, s>r.
\end{equation}
 Indeed, for all~$s\geqslant  r_0,$ since $\|f\|_q\leqslant 1,$ we have
\begin{equation}\label{e1}
\fint_{Q_s}f^qdx = \frac{1}{s^n}\int_{Q_s}f^qdx\leqslant\frac{1}{r_0^n}\int_{\mathbb{R}^n}f^qdx
\leqslant  \frac{1}{r_0^n}=\lambda_0^q\leqslant \lambda^q.
\end{equation}
On the other hand, by Lebesgue's differentiation theorem, we have
\begin{equation*}
\lim_{s\to 0}\fint_{Q_s}f^qdx=f^q(x_0)>\lambda^q.
\end{equation*}
Since $\fint_{Q_s}f^qdx$ is continuous with respect to~$s$, by the intermediate value theorem,  there exists~$0<\overline{r}\leqslant r_0$~such that
$
\fint_{Q_{\overline{r}}}f^qdx=\lambda^q.
$
Set
\[
r=r(x_0)=\sup\{\overline{r}\in (0,r_0]:\fint_{Q_{\overline{r}}}f^qdx=\lambda^q\}.
\]
 Then $0<r\leqslant r_0$ and
\begin{equation}\label{adde1}
\fint_{Q_{ {r}}}f^qdx=\lambda^q.
\end{equation}
To  prove (\ref{e}), suppose on the contrary that  there exists~$s_0>r$~such that
\begin{equation*}
\fint_{Q_{s_{0}}}f^qdx>\fint_{Q_{r}}f^qdx=\lambda^q.
\end{equation*}
Then  by~(\ref{e1}) we must have $s_0<r_0$.
So, again by the intermediate value theorem, there would exist  $r_2\in (s_0, r_0]$~such that
$\fint_{Q_{r_2}}f^qdx=\lambda^q,$ which implies $r_2\leqslant r<s_0,$  a desired contradiction.  \\

{\bf Step 2.}
By (\ref{e}), we obtain
\begin{equation}\label{e2}
\begin{split}
\lambda^q=\fint_{Q_r}f^qdx\leqslant\frac{1}{r^n}\int_{\mathbb{R}^n}f^qdx \leqslant \frac{1}{r^n},\quad r\leqslant\lambda^{-q/n},\\
 \omega(r)\leqslant \omega(\lambda^{-q/n}), \quad \lambda^{\omega(r)}\leqslant\lambda^{\omega(\lambda^{-q/n})}= {a}(\lambda)
 \end{split}
\end{equation}
and thus
\begin{equation}\label{e30}
\begin{split}
 \left(\fint_{Q_{lr}}fdx\right)^{q+\omega(r)} &= \left(\fint_{Q_{lr}}fdx\right)^{q+\omega(r)-1}\fint_{Q_{lr}}fdx \\&
  \leqslant\left(\fint_{Q_{lr}}f^qdx\right)^{\frac{q+\omega(r)-1}{q}}\fint_{Q_{lr}}fdx  \\&
 \leqslant \lambda^{q+\omega(r)-1}\fint_{Q_{lr}}fdx
 \leqslant  {a}(\lambda)\lambda^{q-1} \fint_{Q_{lr}}fdx.
\end{split}
\end{equation}
 Note that
\begin{equation*}
\fint_{Q_{lr}}fdx\leqslant \eta+\frac{1}{|Q_{lr}|}\int_{Q_{lr}\cap E(\eta )}fdx \qquad \forall\, \eta>0,
\end{equation*}
and let
\begin{equation}\label{def-eta}
\eta=\eta(\lambda)=\frac{\lambda}{4c_0a(\lambda)}=\frac{s(\lambda)}{4c_0}.
\end{equation}
Then from~(\ref{e30}), it follows that
\begin{equation}\label{estif0}
c_0\left(\fint_{Q_{lr}}fdx\right)^{q+\omega(r)}\leqslant\frac{\lambda^q}{4}+\frac{c_0{a}(\lambda)\lambda^{q-1}}{|Q_{lr}|}\int_{Q_{lr}\cap E(\eta)}fdx.
\end{equation}
If we set $D(\lambda)=\{x\in\mathbb{R}^n : g(x)>\lambda\},$~then, with $\epsilon=\frac{\lambda}{(4c_1)^{1/q}}$, we have
\begin{equation}\label{estig0}
c_1\fint_{Q_{lr}}g^qdx\leqslant  c_1\epsilon^q  +\frac{c_1}{|Q_{lr}|}\int_{Q_{lr}\cap D(\epsilon)}g^qdx=\frac {\lambda^q}{4}   +\frac{c_1}{|Q_{lr}|}\int_{Q_{lr}\cap D(\epsilon)}g^qdx.
\end{equation}
Now, combining (\ref{gecod0}), (\ref{estif0})~and~(\ref{estig0}),~we  obtain
\begin{align*}
\lambda^q=\fint_{Q_r}f^qdx\leqslant \frac{\lambda^q}{2}+\frac{c_0 {a}(\lambda)\lambda^{q-1}}{ |Q_{lr}|}\int_{Q_{lr}\cap E(\eta)}fdx+\frac{c_1}{|Q_{lr}|}\int_{Q_{lr}\cap D( \epsilon)}g^qdx,
\end{align*}
which gives
\begin{align*}
 \fint_{Q_r}f^qdx\leqslant  \frac{2c_0 {a}(\lambda)\lambda^{q-1}}{|Q_{lr}|}\int_{Q_{lr}\cap E(\eta)}fdx+\frac{2c_1}{|Q_{lr}|}\int_{Q_{lr}\cap D(\epsilon)}g^qdx.
\end{align*}
Hence, using \eqref{e} with $s=4lr$, we have
\begin{align*}
\fint_{Q_{4lr}}f^qdx\leqslant\fint_{Q_r}f^qdx\leqslant  \frac{2c_0 {a}(\lambda)\lambda^{q-1}}{|Q_{lr}|}\int_{Q_{lr}\cap E(\eta)}fdx+\frac{2c_1}{|Q_{lr}|}\int_{Q_{lr}\cap D(\epsilon)}g^qdx,
\end{align*}
which  implies that
\begin{equation}\label{esti50}
\int_{Q_{4lr}}f^qdx\leqslant 2 c_0 4^n {a}(\lambda)\lambda^{q-1}\int_{Q_{lr}\cap E(\eta)}fdx+2c_1 4^{n}\int_{Q_{lr}\cap D(\epsilon)}g^qdx.
\end{equation}
We have thus  proved that for almost all~$x_0\in E(\lambda),$~there exists a cube $Q_r(x_0)$ with  $r=r(x_0)\in (0,r_0]$~for which~(\ref{esti50})~holds. By Vitali's covering theorem there exists a countable family of disjoint cubes~$Q_{lr_j}$~such that
\begin{equation*}
E(\lambda)\subset\bigcup Q_{4lr_j}.
\end{equation*}
Hence by~(\ref{esti50})~we have that, for all $\lambda\geqslant \lambda_0$,
\begin{equation}\label{estlev0}
\int_{E(\lambda)}f^qdx\leqslant  2 c_0 4^n {a}(\lambda)\lambda^{q-1} \int_{E(\eta)}fdx+2c_1 4^{n} \int_{D(\epsilon)}g^qdx.
\end{equation}

{\bf Step 3.} We first observe that if $h\in L^1(\mathbb R^n)$ and $\rho\in C^1([\lambda_0,\infty))$ then
\begin{equation}\label{intequi10}
\int_{E(\lambda_0)}h\rho(f)dx= \rho(\lambda_0)\int_{E(\lambda_0)}h dx+\int_{\lambda_0}^\infty\rho^\prime(\lambda)
\int_{E(\lambda)}h dxd\lambda.
\end{equation}
Indeed, this follows by Fubini's theorem because
\[
\begin{split}\int_{E(\lambda_0)}&h\rho(f)dx = \int_{E(\lambda_0)}h\left(\int_{\lambda_0}^f\rho^\prime(\lambda)d\lambda+\rho(\lambda_0)\right)dx = \rho(\lambda_0)\int_{E(\lambda_0)}h dx+\int_{E(\lambda_0)}h\int_{\lambda_0}^f\rho^\prime(\lambda)d\lambda dx\\
&= \rho(\lambda_0)\int_{E(\lambda_0)}h dx+\iint_{E }h(x) \rho^\prime(\lambda)d\lambda dx = \rho(\lambda_0)\int_{E(\lambda_0)}h dx+\int_{\lambda_0}^\infty\rho^\prime(\lambda)
\int_{E(\lambda)} hdxd\lambda, \end{split}
\]
where  $ E=
\{(x,\lambda): x\in E(\lambda_0), \, \lambda_0<\lambda<f(x)\} =\{(x,\lambda): \lambda_0<\lambda<+\infty,\, x\in E(\lambda)\}.$
We now apply
(\ref{estlev0})~and~(\ref{intequi10}) with  $h=f^q$ and $\varphi=\varphi_{\theta,c_0}(\lambda)$, noticing $\varphi'(\lambda)\geqslant 0$,   to have
\begin{align}\label{esti6}
\int_{E(\lambda_0)}f^q {\varphi}(f)dx&\leqslant 2c_0 4^n  \int_{\lambda_0}^{\infty}
 {\varphi}^\prime(\lambda) {a}(\lambda)\lambda^{q-1}\int_{E(\eta(\lambda))} fdxd\lambda
\nonumber\\
&~~~~+2c_1 4^n
\int_{\lambda_0}^{\infty} {\varphi}^\prime(\lambda)\int_{D(\epsilon)}g^qdxd\lambda+ {\varphi}(\lambda_0)
\int_{\mathbb{R}^n}f^qdx,
\end{align}
where $\eta(\lambda)=\frac{s(\lambda)}{4c_0}$ and $\epsilon=\frac{\lambda}{(4c_1)^{1/q}}.$
Using the substitution $\lambda=\lambda(4c_0\eta),$ where $\lambda(s)$ is the inverse function of $s(\lambda)$ defined above,  by (\ref{intequi10}) with $\rho=\psi,$ the first integral on the right-hand side of (\ref{esti6}) can be written as
\[
\int_{\lambda_0}^{\infty}
 {\varphi}^\prime(\lambda) {a}(\lambda)\lambda^{q-1}\int_{E(\eta(\lambda))} fdxd\lambda=\int_{\eta_0}^{\infty} \psi'(\eta) \int_{E(\eta)} fdxd\eta=\int_{E(\eta_0)}\psi(f)f  dx,
\]
where $\eta_0=\frac{s_0}{4c_0}$ and
 \[
 \psi(\eta)= \int_{4c_0\eta_0}^{4c_0\eta}
 {\varphi}^\prime(\lambda(s) ) {a}(\lambda(s) ) \lambda^{q-1}(s)  \lambda'(s) ds \; \quad\forall\eta\geqslant \eta_0.
 \]
 By (\ref{d-lambda}) and (\ref{phidef0}), we have
 \[
0<\lambda^\prime(s) \leqslant\frac{ {a}(\lambda(s))}{1-\omega_0},\quad   \varphi'(\lambda(s))= \frac{\theta (1-\omega_0)\varphi(\lambda(s))}{\lambda a^q(\lambda(s))} \quad\forall \, s\geqslant s_0.
 \]
Also, by  the relation $\lambda(s)=s {a}(\lambda(s))$,  we have
\begin{equation}\label{addpsi1}
 {\varphi}^\prime(\lambda(s) ) {a}(\lambda(s) ) \lambda^{q-1}(s)  \lambda'(s)  \leqslant  \theta  {\varphi}( \lambda(s))s^{q-2}.
\end{equation}
Thus
\begin{align}\label{addpsi2}
 {\psi}(\eta)\leqslant \theta  {\varphi}( \lambda(4c_0\eta))  \int_{4c_0\eta_0}^{4c_0\eta} s^{q-2}\,ds \leqslant   \frac{\theta (4c_0)^{q-1}{\varphi}( \lambda(4c_0\eta))   }{q-1} \eta^{q-1} \quad\forall\, \eta\geqslant \eta_{0}.
\end{align}
On the one hand, by (\ref{ass-1}), we have
\[
\frac{{\varphi}( \lambda(4c_0\eta)) }{\varphi(\eta)}= \exp\left (\theta(1-\omega_0)\int_\eta^{\lambda(4c_0\eta)} \frac{dt}{ta^q(t)}\right ) \le e^{\theta M(4c_0)}\quad\forall\, \eta\geqslant \eta_{0}=\frac{s_0}{4c_0},
\]
which, by (\ref{addpsi2}), shows that
\[
\psi(\eta)\eta\leqslant   \frac{\theta (4c_0)^{q-1}e^{\theta M(4c_0)} }{q-1}{\varphi}( \eta)   \eta^{q } \quad\forall\, \eta\geqslant \eta_{0}.
\]
This shows that
\begin{equation}\label{addJ1}
\begin{split}
 \int_{\lambda_0}^{\infty}
 {\varphi}^\prime(\lambda) &{a}(\lambda)\lambda^{q-1} \int_{E(\eta(\lambda))} fdxd\lambda   =\int_{E(\eta_0)}\psi(f)f  dx\\& \leqslant  \frac{\theta (4c_0)^{q-1}e^{\theta M(4c_0)} }{q-1} \int_{E(\eta_0)}f^q {\varphi}(f)dx.
\end{split}
\end{equation}
Similarly, by the substitution $t= \epsilon=\frac{\lambda}{(4c_1)^{1/q}},$ we obtain
\begin{equation}\label{esti70}\begin{split}
 \int_{\lambda_0}^{\infty} {\varphi}^\prime(\lambda)\int_{D(\epsilon)}g^qdxd\lambda
&=
 \int_{\frac{\lambda_0}{(4c_1)^{1/q}}}^\infty (4c_1)^{1/q} {\varphi}^\prime((4c_1)^{1/q}t)\int_{D(t)}g^qdxdt
\\
& \leqslant \int_{D(\frac{\lambda_0}{(4c_1)^{1/q}})}g^q {\varphi}((4c_1)^{1/q}g)dx \leqslant \int_{\mathbb{R}^n}g^q {\varphi}((4c_1)^{1/q}g)dx.\end{split}
\end{equation}
Substituting \eqref{addJ1} and \eqref{esti70} into \eqref{esti6} yields
\begin{align}\label{addJ}
\int_{E(\lambda_0)}f^q {\varphi}(f) dx
&\leqslant \frac{4^{n}\theta (4c_0)^{q}e^{\theta M(4c_0)} }{q-1}\int_{E(\eta_0)}f^q {\varphi}(f) dx
+ c_14^{n+1}\int_{\mathbb{R}^n}g^q {\varphi}((4c_1)^{1/q}g) dx + {\varphi}(\lambda_{0})
\nonumber\\
&\leqslant \frac{4^{n}\theta (4c_0)^{q}e^{\theta M(4c_0)} }{q-1}\int_{\mathbb{R}^n}f^q {\varphi}(f) dx
+c_14^{n+1}\int_{\mathbb{R}^n}g^q {\varphi}((4c_1)^{1/q}g)dx+ {\varphi}(\lambda_{0})\nonumber\\
&\leqslant
\frac{4^{n}\theta (4c_0)^{q}e^{\theta M(4c_0)} }{q-1}\int_{E(\lambda_{0})}f^q {\varphi}(f)dx+
\frac{4^{n}\theta (4c_0)^{q}e^{\theta M(4c_0)} }{q-1} \int_{\mathbb{R}^n\backslash E(\lambda_0)}f^q {\varphi}(f)dx
\nonumber\\
&~~~~~+c_14^{n+1}\int_{\mathbb{R}^n}g^q {\varphi}((4c_1)^{1/q}g)dx+ {\varphi}(\lambda_{0})\nonumber\\
&\leqslant
\frac{4^{n}\theta (4c_0)^{q}e^{\theta M(4c_0)} }{q-1} \int_{E(\lambda_{0})}f^q {\varphi}(f)dx+
\bigg(\frac{4^{n}\theta (4c_0)^{q}e^{\theta M(4c_0)} }{q-1}+1\bigg) {\varphi}(\lambda_{0})
\nonumber\\
&~~~~~+c_14^{n+1}\int_{\mathbb{R}^n}g^q {\varphi}((4c_1)^{1/q}g)dx.
\end{align}
Since  $\frac{4^{n}\theta (4c_0)^{q}e^{\theta M(4c_0)} }{q-1}=\frac{1}{2},
$
 it follows   from \eqref{addJ} that
\begin{align*}
\int_{E(\lambda_0)}f^q {\varphi}(f)dx
\leqslant 2c_1 4^{n+1}\int_{\mathbb{R}^n}g^q {\varphi}((4c_1)^{1/q}g)dx+3 {\varphi}(\lambda_{0}).
\end{align*}
Consequently, we obtain
\begin{equation}\label{step3}
\begin{split}
\int_{\mathbb R^n} f^q\varphi(f)\,dx &=\int_{E(\lambda_0)}f^q {\varphi}(f)dx+\int_{\mathbb R^n\setminus E(\lambda_0)}f^q {\varphi}(f)dx\\&\leqslant \int_{E(\lambda_0)}f^q {\varphi}(f)dx+\varphi(\lambda_0) \int_{\mathbb R^n\setminus E(\lambda_0)}f^q dx \\
& \leqslant 2c_1 4^{n+1}\int_{\mathbb{R}^n}g^q {\varphi}((4c_1)^{1/q}g)dx+4 {\varphi}(\lambda_{0}).
 \end{split}
\end{equation}

{\bf Step 4.} Finally, by  (\ref{step4}) and (\ref{step3}), it follows that
\begin{equation}\label{step5}
\begin{split}
\int_{\mathbb R^n} f^q\varphi(f)\,dx &  \leqslant 2c_1 4^{n+1}\int_{\mathbb{R}^n}g^q {\varphi}((4c_1)^{1/q}g)dx+4 {\varphi}(\lambda_{0})\\
&  \leqslant 2c_1 4^{n+1}L((4c_1)^{1/q})\int_{\mathbb{R}^n}g^q {\varphi}(g)dx+4 {\varphi}(\lambda_{0}),
 \end{split}
\end{equation}
which proves  (\ref{result0}) with constant $C=2c_1 4^{n+1}L((4c_1)^{1/q}).$
\end{proof}

\begin{remark}\label{remark-21} Assume, instead, that $\|f\|_q>1$ and (\ref{gecod0}) is satisfied. Then   the functions $\tilde f=f/\|f\|_q$ and $\tilde g=g/\|f\|_q$ will satisfy (\ref{gecod0}) with constant $c_0$ replaced by the new constant $\tilde{c}_0=c_0\|f\|_q.$ With $\varphi=\varphi_{\theta,\tilde{c}_0}$, where $\theta$ is determined by (\ref{theco0}) with $c_0=\tilde{c}_0$, we have
\[
\int_{\mathbb{R}^n}f^q {\varphi}(\tilde f)dx\leqslant 4 {\varphi}(\lambda_0)\|f\|_q^q +C\int_{\mathbb{R}^n}g^q {\varphi}(\tilde g)dx.
\]
From this, by (\ref{step4}) and (\ref{step2}), we establish  the general high integrability estimate:
\begin{equation}\label{result03}
\int_{\mathbb{R}^n}f^q \Phi(f)dx\leqslant C_1\int_{\mathbb{R}^n}f^q dx  +C_2 \int_{\mathbb{R}^n}g^q  \Phi(g)dx,
\end{equation}
where $\Phi(s)$ is a function defined in (\ref{Phi}) with the same $\theta$ determined by (\ref{theco0}) with $c_0=\tilde{c}_0$,  and   $C_1$ and $C_2$  are constants that depend  on  $\|f\|_q.$

Note that the general high integrability estimate  (\ref{result03}) is valid  for all nonnegative $f,g\in L^q(\mathbb R^n)$ as long as  (\ref{gecod0}) is satisfied.
\end{remark}

\section{Sobolev-Poincar\'e inequalities with the variable exponent}

In this section, we present  the modified Sobolev-Poincar\'e inequalities for functions  with the variable exponent. To this end, we first prove some uniform  Sobolev-Poincar\'e inequalities in the ball~$B_r=B_r(x_0)\subset \mathbb R^n.$

\begin{lemma} Let $\beta\geqslant 1.$   Then for all $r>0, \, s\in [1,\beta]$ and $t\in [1,\frac{n}{n-1}],$  the inequality
\begin{equation}\label{poi}
\fint_{B_r}\left|\frac{u-\lambda}{r}\right|^{st}dx\leqslant C_0\left(\fint_{B_r}|\nabla  u|^sdx\right)^t
\end{equation}
holds in the following two cases:
\begin{itemize}
\item[\rm (i)]  For all $ u\in W^{1,s}(B_r)$ with $\lambda=\fint_{B_r}udx,$ where $C_0=c(n,\beta)>0$ is a constant.
\item[\rm (ii)] For  $\lambda=0$ and  all $u\in W^{1,s}(B_r)$ satisfying  $|\{x\in B_r:u(x)=0\}|\geqslant  \gamma |B_r|$ for a given constant $\gamma\in (0,1),$ where $C_0=2c(n,\beta)/\gamma $ with $c(n,\beta)$ being the constant in {\rm Case  (i).}
\end{itemize}
\end{lemma}

\begin{proof} By rescaling $\tilde u(y)=u(ry),$ it suffices to prove   (\ref{poi}) for $r=1$ in both Cases (i) and (ii).  Moreover, by H\"older's inequality: for $t\in [1,\frac{n}{n-1}],$
\[
\left(\fint_{B_1} | u-\lambda|^{st}dx\right)^{\frac{1}{st}} \leqslant \left(\fint_{B_1} | u-\lambda|^{\frac{ns}{n-1}}dx\right)^{\frac{n-1}{ns}},
\]
we only need to prove
\begin{equation}\label{poi-1}
\left(\fint_{B_1} | u-\lambda|^{\frac{ns}{n-1}}dx\right)^{\frac{n-1}{ns}} \leqslant C_0 \left(\fint_{B_1} |\nabla u|^s dx\right)^{\frac{1}{s}}.
\end{equation}

{\bf Step 1:} We prove  (\ref{poi-1})  for all $u\in W^{1,s}(B_1)$  with $\lambda=(u)_{B_1}=\fint_{B_1}udx.$  Let $q$ be such that $\frac{ns}{n-1}=\frac{nq}{n-q}=q^*,$ that is, $q=\frac{ns}{n+s-1}. $  Then $1\leqslant q\leqslant \frac{n\beta}{n+\beta-1}<n$ and $q\leqslant s.$  We claim the uniform  Sobolev inequality: for all $1\leqslant q\leqslant \frac{n\beta}{n+\beta-1}$ and $w\in W^{1,q}(B_1),$
 \begin{equation}\label{emb69}
 \|w\|_{L^{q^*}(B_1)} \leqslant C_1(n,\beta)  \|w\|_{W^{1,q}(B_1)}.
 \end{equation}
We present a direct proof of this inequality using the special feature of unit ball $B_1.$ First assume $w\in C^1(\bar B_1)$ and define
\[
\tilde w(x)=\begin{cases} w(x) & |x|\leqslant 1,\\
w(\frac{x}{|x|^2}) & 1\leqslant |x|\leqslant 2.
\end{cases}
\]
Then $\tilde w\in W^{1,\infty}(B_2).$ Note that the function $y=\psi(x)=\frac{x}{|x|^2}$ has the inverse $x= \psi(y)$ with $\nabla \psi(x)=|x|^{-2} (I-2\frac{x}{|x|}\otimes \frac{x}{|x|})$ and $\det \nabla\psi(x)=-|x|^{-2n}.$ Thus  $|\nabla \tilde w(x)|=|x|^{-2}|\nabla w(\frac{x}{|x|^2})|$ and $dx=|\det\nabla\psi(y)|dy=|y|^{-2n}dy$ on $1<|x|<2.$ Hence it can be shown that
 \begin{equation}\label{emb70}
\|\tilde w\|_{L^q(B_2)}\le (1+4^n)\| w\|_{L^q(B_1)},\quad\|\nabla \tilde w\|_{L^q(B_2)}\le (1+4^n)\| \nabla w\|_{L^q(B_1)}.
\end{equation}
Let $\zeta\in W^{1,\infty}_0(B_2)$ be such that $0\leqslant \zeta \leqslant  1,\, |\nabla\zeta |\leqslant  1$ on $B_2$ and $\zeta |_{B_1}\equiv 1.$ Set $v=\zeta \tilde w\in W_0^{1,\infty}(B_2).$  We apply the Gagliardo-Nirenberg-Sobolev-Poincar\'e  inequality:
\[
 \|v\|_{L^{q^*}(B_2)} \leqslant  \frac{q(n-1)}{n-q} \| \nabla v\|_{L^{q}(B_2)}
\]
and use $\nabla v=\zeta\nabla\tilde w+ \tilde w \nabla\zeta $ and (\ref{emb70}) to obtain
 \[
 \|w\|_{L^{q^*}(B_1)}  \leqslant  \|v\|_{L^{q^*}(B_2)}\leqslant  \frac{q(n-1)}{n-q} \| \nabla v\|_{L^{q}(B_2)} \]\[ \leqslant (1+4^n) \frac{q(n-1)}{n-q} \|w\|_{W^{1,q}(B_1)}  \leqslant (1+4^n) (n+\beta-1) \|w\|_{W^{1,q}(B_1)},
 \]
which proves (\ref{emb69}) for smooth $w$ with constant $C_1(n,\beta) =(1+4^n) (n+\beta-1).$ The general case follows by the standard density argument.  We next claim the uniform Sobolev-Poincar\'e inequality: for all $1\leqslant q\leqslant \frac{n\beta}{n+\beta-1}$ and $u\in W^{1,q}(B_1),$
\begin{equation}\label{emb72}
\left(\fint_{B_1} | u-\lambda|^{q^*} \right)^{\frac{1}{q^*}}  \leqslant c(n,\beta) \left(\fint_{B_1} |\nabla u|^q  \right)^{\frac{1}{q}},
\end{equation}
from which we deduce (\ref{poi-1})  as follows:
 \[
 \left(\fint_{B_1} | u -\lambda|^{\frac{ns}{n-1}} \right)^{\frac{n-1}{ns}}  \leqslant c(n,\beta) \left(\fint_{B_1} |\nabla u|^q  \right)^{\frac{1}{q}} \leqslant c(n,\beta) \left(\fint_{B_1} |\nabla u|^s  \right)^{\frac{1}{s}}.
\]
To prove (\ref{emb72}),  suppose, on the contrary, that there exist $1\leqslant q_j\leqslant \frac{n\beta}{n+\beta-1}$ and $u_j\in W^{1,q_j}(B_1)$ such that
\begin{equation}\label{emb73}
\left(\fint_{B_1} | u_j-\lambda_j|^{q_j^*} \right)^{\frac{1}{q_j^*}} =1,\quad  \lim_{j\to\infty} \left(\fint_{B_1} |\nabla u_j|^{q_j}  \right)^{\frac{1}{q_j}}=0.
\end{equation}
We  assume $q_j\to\bar q\in [1,\frac{n\beta}{n+\beta-1}].$ Then it is possible to select  numbers $q$ and $r$ such that
\[
1\leqslant q\leqslant q_j\leqslant r<q^*\leqslant q_j^* \qquad \forall\,j>>1.
\]
Then $\{u_j-\lambda_j\}$ is a bounded sequence in $W^{1,q}(B_1).$ By the Sobolev compact embedding theorem, we assume (possibly by taking a subsequence) that $u_j-\lambda_j\to \bar u$ in $L^r(B_1).$ Since $\nabla (u_j-\lambda_j)=\nabla u_j \to 0$ in $L^q(B_1)$, it follows that $\nabla \bar u =0$ on $B_1$; thus, as $\fint_{B_1} \bar u dx=0$, we must have $\bar u=0$ on $B_1.$ Hence
\[
\lim_{j\to\infty}
\left(\fint_{B_1} |u_j-\lambda_j|^r dx\right)^{1/r} =0.
\]
 On the other hand, applying (\ref{emb69}) with $q=q_j$ and $w=u_j-\lambda_j$, thanks to (\ref{emb73}), it follows that for some constant $\epsilon_0>0$,
\[
\|u_j-\lambda_j\|_{L^{q_j}(B_1)} \geqslant \epsilon_0  \qquad \forall\,j>>1,
\]
which yields  a desired contradiction:
\[
\left(\fint_{B_1} |u_j-\lambda_j|^r dx\right)^{1/r}\geqslant \left(\fint_{B_1} |u_j-\lambda_j|^{q_j} dx\right)^{1/q_j}\geqslant \epsilon_0'>0 \qquad \forall\,j>>1.
\]

{\bf Step 2:} We prove (\ref{poi-1}) with $\lambda=0$ for all  $
 u\in W^{1,s}(B_1)$ satisfying
 \[
 |\{x\in B_1:u(x)=0\}|\geqslant  \gamma |B_1|,
 \]
   where $\gamma\in (0,1)$ is a constant.  Given such a function $u$, let $E=\{x\in B_1:u(x)=0\}.$ From the inequality (\ref{poi-1})  just proved in Step 1, we have
  \[
\begin{split}
|E|^{\frac{n-1}{ns}} |(u)_{B_1}|&=\left(\int_{E} |(u)_{B_1}|^{\frac{ns}{n-1}}dx\right)^{\frac{n-1}{ns}}  =\left(\int_{E} | u -(u)_{B_1}|^{\frac{ns}{n-1}}dx\right)^{\frac{n-1}{ns}}\\
& \leqslant \left(\int_{B_1} | u -(u)_{B_1}|^{\frac{ns}{n-1}}dx\right)^{\frac{n-1}{ns}}   \leqslant  c(n,\beta) |B_1|^{\frac{n-1}{ns}} \left(\fint_{B_1} |\nabla u|^s dx\right)^{\frac{1}{s}}.\end{split}
\]
This, together with $|E|\geqslant \gamma|B_1|,$ yields  that
\[
|(u)_{B_1}| \leqslant c(n,\beta)\Big (\frac{|B_1|}{|E|}\Big)^{\frac{n-1}{ns}} \left(\fint_{B_1} |\nabla u|^s dx\right)^{\frac{1}{s}}\leqslant  \frac{c(n,\beta)}{\gamma}  \left(\fint_{B_1} |\nabla u|^s dx\right)^{\frac{1}{s}}.
\]
Thus
\[
\|(u)_{B_1}\|_{L^{\frac{ns}{n-1}}(B_1)} =|B_1|^{\frac{n-1}{ns}}|(u)_{B_1}| \leqslant \frac{c(n,\beta)}{\gamma}  |B_1|^{\frac{n-1}{ns}} \left(\fint_{B_1} |\nabla u|^s dx\right)^{\frac{1}{s}}.
\]
Finally, we obtain
\[
\begin{split}
 \|u\|_{L^{\frac{ns}{n-1}}(B_1)} &\leqslant  \|u-(u)_{B_1}\|_{L^{\frac{ns}{n-1}}(B_1)}  +  \|(u)_{B_1}\|_{L^{\frac{ns}{n-1}}(B_1)}    \\
 & \leqslant \Big [ c(n,\beta)+\frac{c(n,\beta)}{\gamma}\Big ]   |B_1|^{\frac{n-1}{ns}}  \left(\fint_{B_1} |\nabla u|^s dx\right)^{\frac{1}{s}} \\ &\leqslant  \frac{2c(n,\beta)}{\gamma}  |B_1|^{\frac{n-1}{ns}}  \left(\fint_{B_1} |\nabla u|^s dx\right)^{\frac{1}{s}},
 \end{split}
\]
 which  proves (\ref{poi-1})  with $\lambda=0$ and $C_0=\frac{2c(n,\beta)}{\gamma}.$
\end{proof}

\begin{theorem}\label{gesobpoi1} Let $p(x)$ satisfy (\ref{omeg}) with $\alpha>1$ and let $r_0\in (0,1)$ further satisfy  $\omega(r_0)< \frac{1}{n-1}.$ Then for all
\begin{equation}\label{cond-q}
0<r\leqslant r_0, \quad 1\leqslant q \leqslant \min\Big\{\alpha,\, \frac{n}{n-1} -\omega(r_0)\Big \},
\end{equation}
 the inequality
\begin{equation}\label{poin1}
\fint_{B_r }\left|\frac{u-\lambda}{r}\right|^{p(x)}dx\leqslant c+c\left(\fint_{B_r }|\nabla  u|^{\frac{p(x)}{q}}dx\right)^{q+{\omega(r)}}
\end{equation}
holds for all $u\in W^{1,p(x)}(B_r)$ and $\lambda=(u)_{B_r},$ where $c=c(n, \beta)>0$ is a constant. Moreover, the same inequality holds with $\lambda=0$ for all $u\in W^{1,s}(B_r)$ satisfying  $|\{x\in B_r:u(x)=0\}|\geqslant  \gamma |B_r|$ for a given constant $\gamma\in (0,1),$ where $c=c(n,\beta, \gamma)>0$  is a constant.
\end{theorem}
\begin{proof}
Let $p_1=\min\limits_{\bar{B}_r}p,\,  p_2=\max\limits_{\bar{B}_r}p,\, s=\frac{p_1}{q}$ and $t=\frac{qp_2}{p_1}.$ Clearly,
\[
\alpha\leqslant p_1\leqslant p(x)\leqslant p_2\leqslant \beta \quad\mbox{and}\quad  p_2\leqslant p_1+\omega(r).
\]
Thus, by (\ref{cond-q}), we have
\[
1\leqslant s\leqslant \beta, \quad 1\leqslant t\leqslant \frac{p_1+\omega(r)}{p_1}q\leqslant q+ \omega(r) \leqslant q+ \omega(r_0)   \leqslant\frac{n}{n-1}.
\]
From (\ref{poi}) and since $(1+a)^k\leqslant 2^k(1+a^k)$ for all nonnegative $a$ and $k$,  we obtain
\begin{eqnarray*}
\fint_{B_r}\left|\frac{u-\lambda}{r}\right|^{p(x)}dx&\leqslant&1+\fint_{B_r}\left|\frac{u-\lambda}{r}\right|^{p_2}dx\\
&\leqslant&1+C_0\left(\fint_{B_r}|\nabla  u|^{\frac{p_1}{q}}dx\right)^t\\
&\leqslant&1+C_0\left(1+\fint_{B_r}|\nabla  u|^{\frac{p(x)}{q}}dx\right)^{q+ \omega(r)} \\
&\leqslant&1+2^{\frac{n}{n-1}}C_0\left[1+\left(\fint_{B_r}|\nabla  u|^{\frac{p(x)}{q}}dx\right)^{q+ \omega(r)} \right ],
\end{eqnarray*}
which proves (\ref{poin1}) with $c=1+2^{\frac{n}{n-1}}C_0.$
\end{proof}

\section{Proof  of  the main result}

In this section, we assume that $p(x)$ satisfies the assumptions (\ref{omeg}) and (\ref{ass-0}) for some $\alpha>2.$ Denote
\[
p'(x)=\frac{p(x)}{p(x)-1}.
\]

\subsection{Cacciopolli's inequality}\label{cacc}

  \begin{theorem}  \label{cacci}  Let    $u\in W_0^{1,p(x)}(\Omega)$ be a weak  solution of~(\ref{equ}). {   We extend $\vec{b}, \vec{h}$ and $u$ by zero outside $\Omega.$  Then there exists a constant $c =c(\alpha,\beta)>0$ such that
\begin{equation}\label{cacciineq}
\int_{B_{r/2}(x_0)}|\nabla  u|^{p(x)}dx\leqslant c\int_{B_r(x_0)}\left(\left|\frac{u-\lambda}{r}\right|^{p(x)}+|\vec{b}|^{\frac{p(x)}{p(x)-2}} +|\vec{h}|^{p^\prime(x)}\right)dx,
\end{equation}
for all balls $B_r=B_r(x_0)$ with $0<r\leqslant 1$,  where
\begin{equation} \label{def-lambda}
\lambda=
\begin{cases}
\fint_{B_r}udx,&\text{if}~B_{3r/4}\subset\Omega,\\
0,&\text{otherwise}.
\end{cases}
\end{equation}
}
\end{theorem}

 \begin{proof} Let $B_r=B_r(x_0)$ with $0<r\leqslant 1$ be any ball. Let $\chi\in C^{\infty}_0(B_{3r/4})$ be such that
\begin{equation*}
0\leqslant\chi\leqslant 1,\quad \chi|_{B_{r/2}}\equiv 1,\quad |\nabla \chi|\leqslant 5r^{-1}.
\end{equation*}
Consider function $\varphi=\chi^\beta(u-\lambda),$ where
\begin{equation*}
\lambda=
\begin{cases}
\fint_{B_r}udx,&\text{if}~B_{3r/4}\subset\Omega,\\
0,&\text{otherwise}.
\end{cases}
\end{equation*}
It can be shown that $\varphi\in W^{1,p(x)}_0(\Omega)$ is a legitimate    test function for  (\ref{equ}).
Hence, we have
\begin{equation}\label{esti1}
\begin{split} \int_{B_r}|\nabla  u|^{p(x)}\chi^{\beta}dx  \leqslant  & \, 5\beta\int_{B_r}\chi^{\beta-1}(|\nabla  u|^{p(x)-1}  +|\vec{h}|)\left |\frac{u-\lambda}{r}\right |dx\\
& +\int_{B_r}\chi^\beta \vec{h}\cdot\nabla  udx \\
 &+\int_{B_r}\chi^\beta(u-\lambda) \vec{b}(x)\cdot\nabla  udx \\
   = \, &J_{1}+J_{2}+J_{3}.
 \end{split}
\end{equation}
First, we estimate the terms $J_{1}$ and $J_{2}$. By  Young's inequality with~$\varepsilon\in(0,1)$ we have
\begin{eqnarray}\label{addJ_{11}}
\chi^{\beta-1}|\nabla  u|^{p(x)-1}\left|\frac{u-\lambda}{r}\right|&\leqslant&
\varepsilon^{p^\prime(x)}\chi^{(\beta-1)p^\prime(x)}|\nabla  u|^{p(x)}+\varepsilon^{-p(x)}\left|\frac{u-\lambda}{r}\right|^{p(x)}\nonumber\\
&\leqslant&\varepsilon^{\beta^\prime}\chi^\beta|\nabla  u|^{p(x)}+\varepsilon^{-\beta}\left|\frac{u-\lambda}{r}\right|^{p(x)}.
\end{eqnarray}
Here we have used the facts that $0\leqslant \chi\leqslant 1$, $p(x)\leqslant \beta $ and $(\beta-1)p'(x)\geqslant (\beta-1)\beta^\prime=\beta.$ Similarly, we have
\begin{equation}\label{addJ_{12}}
\begin{split} |\vec{h}\chi^\beta\cdot\nabla  u| &\leqslant\varepsilon^{p(x)}\chi^\beta|\nabla  u|^{p(x)}+\varepsilon^{-p'(x)}\chi^\beta|\vec{h}|^{p^\prime(x)}\\&\leqslant\varepsilon^{\alpha}\chi^\beta|\nabla  u|^{p(x)}+\varepsilon^{-\alpha'} |\vec{h}|^{p^\prime(x)},\end{split}
\end{equation}
and
\begin{equation}\label{addJ_{13}}
 |\vec{h}| \left |\frac{u-\lambda}{r}\right|\leqslant|\vec{h}|^{p^\prime(x)}+\left|\frac{u-\lambda}{r}\right|^{p(x)}.
\end{equation}
Thus, combining \eqref{addJ_{11}}-\eqref{addJ_{13}}, we have
\begin{align}\label{J_{1}}
J_{1}&\leqslant\varepsilon^{\beta^\prime}\int_{B_{r}}\chi^\beta|\nabla  u|^{p(x)}dx+\Big (1+\varepsilon^{-\beta}\Big ) \int_{B_{r}}\left|\frac{u-\lambda}{r}\right|^{p(x)}dx +\int_{B_{r}}|\vec{h}|^{p^\prime(x)}dx, \\
J_{2}&\leqslant\varepsilon^{\alpha} \int_{B_{r}}\chi^\beta|\nabla  u|^{p(x)}dx+\varepsilon^{-\alpha'}\int_{B_{r}} |\vec{h}|^{p^\prime(x)}dx.\nonumber
\end{align}
 We now estimate  the term $J_{3}$. Again, by Young's inequality with $\varepsilon\in(0,1)$ we have
\begin{equation}\label{j3}
\begin{split}
J_{3} &\leqslant\int_{B_r}\varepsilon^{p(x)}|\nabla  u|^{p(x)}\chi^{\beta}dx+\int_{B_r}\varepsilon^{-p^\prime(x)}\chi^{\beta}|\vec{b}(x)(u-\lambda)|^{p^\prime(x)}dx\\
&\leqslant \varepsilon^{\alpha} \int_{B_r}|\nabla  u|^{p(x)}\chi^{\beta}dx+\varepsilon^{-\alpha'}  \int_{B_r}|\vec{b}(x)(u-\lambda)|^{p^\prime(x)}dx.\end{split}
\end{equation}
Thus if $0<r\leqslant 1$ then
\begin{align*}
\int_{B_r}|\vec{b}(x)(u-\lambda)|^{p^\prime(x)}dx&=\int_{B_r}|r\vec{b}(x)|^{p^\prime(x)}|\frac{u-\lambda}{r}|^{p^\prime(x)}dx\\
&\leqslant\int_{B_r}\Big [|r\vec{b}(x)|^{\frac{p(x)}{p(x)-2}}+|\frac{u-\lambda}{r}|^{p(x)} \Big ] dx\\
&\leqslant  \int_{B_r}|\vec{b}(x)|^{\frac{p(x)}{p(x)-2}}dx+\int_{B_r}|\frac{u-\lambda}{r}|^{p(x)}dx.
\end{align*}
 Finally, combining (\ref{J_{1}}) and (\ref{j3}) and choosing a sufficiently small $\varepsilon>0$, we establish (\ref{cacciineq}) with a constant $c=c(\alpha,\beta)>0.$
\end{proof}

\subsection{Reverse H\"older inequality}\label{Hol}

To incorporate  Theorem \ref{cacc}  with constant $\lambda$ defined  by (\ref{def-lambda})  in order to apply  Theorem \ref{gesobpoi1},  we need  some restrictions on the  domains $\Omega\subset\mathbb R^n.$

\begin{definition}\label{cone-Omega} A domain $\Omega\subset\mathbb R^n$ is said to satisfy the  {\rm uniform exterior cone condition} provided that there exist number $\bar  r=\bar{r}(\Omega)>0$ and  cone $C_0=C_0(\Omega)$ with vertex $0$ in $\mathbb R^n$ such that $
B_{\bar{r}}(y)\cap \{ y+RC_0\} \subset \Omega^c$  for all $y\in \partial\Omega,$  where $R=R(y)$ is a rotation.
 \end{definition}

 \begin{lemma} Let $\Omega$ be a bounded domain  satisfying the uniform exterior cone condition with constant $\bar r=\bar{r}(\Omega)$ and cone $C_0=C_0(\Omega)$ as above.    Then there exists a constant $\gamma\in (0,1)$ such that for all $r\in (0,\bar r)$ and $x_0\in\mathbb R^n$ it follows that
\begin{equation} \label{ass-omega}  B_{3r/4}(x_0)\cap \Omega^c\ne \emptyset \Longrightarrow  |B_{r}(x_0)\cap \Omega^c|\geqslant \gamma |B_r|.
\end{equation}
\end{lemma}

  \begin{proof}  Let  $r\in (0,\bar{r})$ and $x_0\in\mathbb R^n$ be such that $B_{3r/4}(x_0)\cap \Omega^c\ne \emptyset.$
If $B_{3r/4}(x_0)\cap\partial\Omega = \emptyset,$  then  $B_{3r/4}(x_0)\subset \Omega^c$ and hence
  \[
|B_{r}(x_0)\cap \Omega^c| \geqslant |B_{3r/4}(x_0)\cap \Omega^c| =|B_{3r/4}(x_0)|= (3/4)^n |B_r|.
\]
If $B_{3r/4}(x_0)\cap\partial\Omega \ne \emptyset$ and let  $y_0\in B_{3r/4}(x_0)\cap\partial\Omega,$ then $B_{r/4}(y_0)\subset B_{r}(x_0)$  and hence
 \[
|B_{r}(x_0)\cap \Omega^c| \geqslant |B_{r/4}(y_0)\cap \Omega^c| \geqslant |B_{r/4}(y_0)\cap (y_0+R_0C_0)|= \epsilon_0 |B_r|,
\]
 where $R_0$ is a rotation and
 \[
  \epsilon_0=\frac{|B_{1/4}(0)\cap  C_0|}{|B_1(0)|}<1/4^n
  \]
   is a positive number depending only on the cone $C_0.$ This proves (\ref{ass-omega})  with $\gamma= \epsilon_0.$
  \end{proof}

We are now in a position to establish the reverse H\"older inequality.

\begin{theorem}\label{reHolder2} Let $\Omega$ be a bounded domain  satisfying the uniform exterior cone condition and  let $ |\vec{b}(x)|^{\frac{p(x)}{p(x)-2}}+ |\vec{h}(x)|^{\frac{p(x)}{p(x)-1}} \in L^1(\Omega).$
Assume  that
  $p(x)$ satisfies   (\ref{omeg}) and (\ref{ass-0}) with $\alpha>2$ and that $r_0\in (0,1)$ further satisfies  $\sqrt{n} r_0 < \bar{r}(\Omega)$ and $\omega(\sqrt{n} r_0)<\frac{1}{n-1}.$   Then there exist constants $c_0$ and $c_1$  such that  the inequality
\begin{equation}\label{reholder2}
\fint_{B_{r/2}(x_0)}|\nabla  u|^{p(x)}dx\leqslant c_0\left(\fint_{B_r(x_0)}|\nabla  u|^{\frac{p(x)}{q}} dx\right)^{q+ \omega(r) }
+c_1\fint_{B_r(x_0)} G(x)dx,
\end{equation}
where $G= 1+ |\vec{b}(x)|^{\frac{p(x)}{p(x)-2}}+ |\vec{h}(x)|^{\frac{p(x)}{p(x)-1}}$, holds for all weak solutions $u\in W^{1,p(x)}_0(\Omega)$ of problem (\ref{equ}),  $r\in (0,r_0)$, $x_0\in\mathbb R^n$ and $1\leqslant q\leqslant  \frac{n}{n-1}- \omega(r_0).$ Here $u$, $\vec{b}$ and $\vec{h}$ are all extended by zero outside $\Omega.$
\end{theorem}

\begin{proof} In Theorem \ref{cacci}~and Theorem \ref{gesobpoi1}, if $B_{3r/4}(x_0)\subset \Omega$ then $\lambda=(u)_{B_r}$; if $B_{3r/4}(x_0)\cap \Omega^c\ne \emptyset$ then, by   (\ref{ass-omega}), we have $|B_{r}(x_0)\cap \Omega^c|\geqslant \gamma |B_r|.$ Finally, in both cases,   Theorem  \ref{cacci} and Theorem \ref{gesobpoi1} are applicable, which proves  (\ref{reholder2}).
\end{proof}

\subsection{Proof of Theorem~\ref{main3}} We extend  $u,\vec{b}$ and $\vec{h}$ to be zero outside $\Omega.$
Take
\begin{equation*}
f=|\nabla  u|^{\frac{p(x)}{q}},\quad G= 1+|\vec{b}(x)|^{\frac{p(x)}{p(x)-2}}+ |\vec{h}(x)|^{\frac{p(x)}{p(x)-1}}, \quad g=G^{\frac{1}{q}}.
\end{equation*}
Let $q$ be any number such that $1<q\leqslant \frac{n}{n-1}-  \omega(r_0),$  where $r_0$ is as specified  in  Theorem~\ref{reHolder2}.  Applying Theorem~\ref{reHolder2} and observing that
\[
Q_r(x_0)\subset B_{{\sqrt{n}\,r}/{2} }(x_0)\subset B_{\sqrt{n} \,r}(x_0) \subset Q_{2\sqrt{n} \, r}(x_0),
\]
 we  obtain (\ref{gecod0}) with $l=2\sqrt{n},$  which, combined with (\ref{result03}), proves  ~(\ref{h-i}) under (\ref{g-i}).
\\
\\

\noindent{\bf Acknowledgments.}
First two authors were supported by NSFC (11301211) and NSF of Jilin Prov. (201500520056JH). Jingya Chen was also supported by Graduate Innovation Fund of Jilin University (2022069).
\\

\noindent{\bf The data availability statement.}
Data sharing not applicable to this article as no datasets were generated
or analyzed during the current study.

\end{document}